\documentclass[12pt]{amsproc}

\title[Classification of Exponential Maps]{A Combinatorial Classification of Postsingularly Finite Complex Exponential Maps}
\author{Bastian Laubner}
\author{ Dierk Schleicher}
\author{ Vlad Vicol}
\address{Department of Mathematics, Malott Hall, Cornell University, Ithaca, NY 14853, USA}
\email{bl232$@$cornell.edu}
\address{School of Engineering and Science, Jacobs University Bremen (formerly International University Bremen), Postfach 750 561, D-28725 Bremen, Germany}
\email{dierk@iu-bremen.de}
\address{Department of Mathematics, University of Southern California,
3620 South Vermont Ave., KAP 108, Los Angeles, California 90089-2532}
\email{vicol@usc.edu}
\dedicatory{Dedicated to the memory of Des Sheiham,\\ our inspiring instructor, valued colleague and great friend}

\usepackage{amscd}
\newtheorem{lemma}{Lemma}[section]
\newtheorem{defn}[lemma]{Definition}
\newtheorem{thm}[lemma]{Theorem}
\newtheorem{prop}[lemma]{Proposition}
\newtheorem{coroll}[lemma]{Corollary}

\newtheorem{algo}[lemma]{Algorithm}

\def\proof{\par\medskip\noindent {\sc Proof. }}
\def\remark{\par\medskip \noindent {\sc Remark. }}
\newcommand{\proofof}[1]{\par\medskip\noindent {{\sc Proof of} #1. }}
\def\qed{\rule{0pt}{3pt}\hfill $\square$}
\def\qedd{\qed}

\newfont{\bbf}{msbm10 at 11pt}

\def\R{\mathbb R}
\def\N{\mathbb N}
\def\C{\mathbb C}
\def\Z{\mathbb Z}
\def\S{\mathbb S}
\def\D{\mathbb D}

\newcommand{\itin}[1]{\mathrm{\underline{#1}}}
\newcommand{\capt}[1]{{\bf(#1)}\\}
\newcommand{\It}[2]{\mathrm{It}(#1 \mid #2)}

\newcommand{\im}{\mathrm{Im}}

\newcommand{\Cbar}{\overline{\C}}
\newcommand{\Cstar}{\C^*}
\newcommand{\disk}{\mathbb {D}}

\newcommand{\Sphere}{{\mathbb S}^2}

\newcommand{\ovl}{\overline}

\renewcommand{\Re}{{\rm Re}}
\renewcommand{\Im}{{\rm Im}}
\newcommand{\El}{E_\lambda}
\newcommand{\Elz}{E_{\lambda_0}}
\newcommand{\Elp}{E_{\lambda'}}
\newcommand{\s}{{\underline s}}
\renewcommand{\t}{\underline t}
\renewcommand{\u}{{\mbox{\tt u}}}

\newcommand{\uu}{\underline \u}

\newfont{\Euler}{eusm10 at 12pt}
\newcommand{\Sym}{ \mathcal{S} }
\newcommand{\sm}{\setminus}
\newcommand{\hide}[1]{}

\newcommand{\ray}{g_\s}
\newcommand{\rayp}{g_{\s'}}
\newcommand{\Ray}{G_\s}

\begin{document}


%%%%%%%%%%%%%%%%%%%%%%%%%%%%%%%%%%%%%%%%%%%%%%%%%%%%%%%%%%%%%%%%%%%%%%%%%%%

\begin{abstract}
We give a combinatorial classification of postsingularly finite
exponential maps in terms of external addresses starting with the
entry $0$. This is an extension of the classification results for
critically preperiodic polynomials \cite{BFH} to exponential maps.
Our proof relies on the topological characterization of
postsingularly finite exponential maps given recently in
\cite{HSS}. Our results illustrate once again the fruitful interplay between combinatorics, topology and complex structure which has often been successful in complex dynamics.
\end{abstract}

\maketitle

%%%%%%%%%%%%%%%%%%%%%%%%%%%%%%%%%%%%%%%%%%%%%%%%%%%%%%%%%%%%%%%%%%%%%%%%%%%
%\newpage
\tableofcontents
%\newpage
%%%%%%%%%%%%%%%%%%%%%%%%%%%%%%%%%%%%%%%%%%%%%%%%%%%%%%%%%%%%%%%%%%%%%%%%%%%%%

\section{Introduction}
\label{sec:intro}

We study the dynamical systems given by iteration of exponential maps $z\mapsto\El(z):=\lambda\exp(z)$ for non-zero complex parameters $\lambda$. The family of exponential maps is the simplest family of transcendental entire functions and has been investigated by many people (see for example \cite{BR,EL,DGH}), often in analogy to quadratic polynomials as
 the simplest family of algebraic entire functions.

The dynamics of iterated holomorphic functions $f\colon\C\to\C$ is determined to a large extent by the dynamics of the singular values: these are values $a\in\C$ which have no neighborhood $U\subset\C$ so that $f$ is an unramified covering over $U$. For polynomials, singular values are critical values. For transcendental functions, singular values can also be asymptotic values, or limit points of critical or asymptotic values. The exponential family is special because it has only one asymptotic value, just like unicritical polynomials (those conjugate to $z\mapsto z^d+c$) have only one critical value.

In any family of iterated holomorphic functions $f_\lambda\colon\C\to\C$, the easiest maps to understand are usually those for which all singular values have finite orbits, i.e.\ the singular orbits are periodic or preperiodic; such maps are called {\em postsingularly finite} (or, for polynomials, postcritically finite). Often they are also the maps which are most important for the structure of parameter space.

The main example is the Mandelbrot set: the {\em Branch Theorem}
\cite{Orsay,Branch} asserts that branch points (in a precise sense)
within the Mandelbrot set are postcritically finite, and the entire
topology of the Mandelbrot set is completely described by them
(under the assumption of local connectivity): if the unique critical
orbit is periodic, the corresponding parameter is the center of a
hyperbolic component, while if the critical orbit is preperiodic,
the parameter is called a ``Misiurewicz point''. For iterated
rational functions, there is a powerful theorem by Thurston
\cite{DH} which helps to understand postcritically finite rational
functions; a variant for polynomials is known as {\em spiders}
\cite{HS}. While Thurston's theorem is deep and powerful, each time
it is applied is usually a theorem in its own right. For instance,
the classification of quadratic polynomials with periodic critical
orbits in \cite{HS}, the classification of general polynomials with
preperiodic critical orbits in \cite{BFH}, and the classification of
general postcritically finite polynomials in \cite{Poirier} are all
derived form Thurston's theorem.

There are many analogies between the bifurcation locus of quadratic
polynomials (the boundary of the Mandelbrot set) and the bifurcation locus of exponential
functions. Again, one expects that much of the structure of the
bifurcation diagram is determined by hyperbolic components and
postsingularly finite exponential maps. For a recent survey about
exponential parameter space, see \cite{CRAS,RS,ExpoPara}.

Since for exponential maps, the singular value $0$ is an omitted value, it can never be periodic, so hyperbolic components have no center. Hyperbolic components have been classified completely in \cite{S}. Postsingularly finite exponential maps thus necessarily have preperiodic singular orbits and are sometimes equivalently called ``postsingularly preperiodic''. Since Thurston's theorem applies only to rational maps, the investigation of postsingularly finite entire functions is much harder. Recently, \cite{HSS} provided an extension of Thur\-ston's theorem specifically to postsingularly finite exponential maps. We make essential use of that theorem. Our main result is a  combinatorial classification of exponential functions $z\mapsto\lambda\exp(z)$ for which the singular value $0$ is preperiodic. Our classification is in terms of preperiodic external addresses, i.e.\ preperiodic sequences over the integers. We should mention that Bergweiler (unpublished) used value distribution theory to estimate the density of postsingularly finite exponential maps.

Our result also contributes to answering (a generalization of) an old question of Euler~\cite{Euler}: for which values $a$ does the limit $a^{a^{a^{\dots}}}$ exist? Euler asked this only for real $a>0$, for which the answer is relatively simple; if $a$ is allowed to be complex, the answer has a very rich structure. In order to be well-defined, the question needs to be rewritten: fixing a branch $\lambda=\log(a)$, then $a^{a^{a^{\dots}}}=e^{\lambda e^{\lambda e^{\dots}}}$, and we are asking for which $\lambda$ the sequence $e^\lambda$, $e^{\lambda e^\lambda}$, $e^{\lambda e^{\lambda e^\lambda}}$, \dots has a limit; except for the final exponentiation step, this is asking for which values of $\lambda$ the exponential map $z\mapsto\lambda e^z$ has a converging singular orbit. The answer to this comes in three parts: (a) convergence in $\C$ without being eventually constant; (b) eventually constant convergence; (c) convergence to $\infty$. Part (a) is easy to answer: this happens iff $\lambda=\mu e^{-\mu}$ with $|\mu|<1$ or $\mu$ a root of unity ($\El$ has an attracting or parabolic fixed point). Part (c) has been answered in \cite{FRS}: the corresponding locus in parameter space consists of uncountably many curves in $\lambda$-space called parameter rays (see Proposition~\ref{prop:pararays}). Finally,  part (b) are exactly the postsingularly finite exponential maps, and their classification is our main result.

\hide{
{\sc Organization of the Paper}. In Section~\ref{sec:defns}, we introduce the necessary background about exponential dynamics; we define dynamic rays and their combinatorial description in terms of external addresses. In Section~\ref{sec:holexpmap2extadd} we associate to every postsingularly finite holomorphic exponential $\El$ map an external address: the singular value is the landing point of a preperiodic dynamic ray, and its external address is a combinatorial coding of $\El$ (there is an ambiguity if several rays land at the singular value). In Section~\ref{sec:combitins} we introduce tools from symbolic dynamics: itineraries and kneading sequences. The main construction is done in Section~\ref{sec:extadd2topexpmap}: to a given external address $\s$, we first define an infinite graph with a continuous self-map. We embed this graph into the two-sphere and show that the graph map can be extended to the entire sphere, yielding a topological exponential map. In Section~\ref{sec:topexp2holexp}, we show that this topological exponential map is related to a unique postsingularly finite holomorphic exponential map in such a way that the singular value is the landing point of the dynamic ray at external address $\s$. We also discuss under which conditions different external addresses describe the same postsingularly finite exponential map. Finally, in Section~\ref{Sec:Parameter}, we show that our classification is natural from the point of view of exponential parameter space and landing properties of parameter rays.
}

This paper grew out of the Bachelor's theses of Bastian and Vlad at International University Bremen in spring 2005. We would like to thank Nikita Selinger and an anonymous referee for many helpful comments.

\section{Definitions and Classification
Theorem}\label{sec:defns}
%Diagram, basic definitions and main statements. + Vlad:

In this section, we introduce the necessary background from exponential and symbolic dynamics, we state our main theorems, and we present a global overview of the argument and thus of the entire paper.

{\sc Notation.}
We set  $\C^* \colon= \C \setminus \{0\}$,\ $\C^{'} := \C^{*}
\setminus \R^{-}$,\ $\Cbar := \C \cup \{
\infty\}$ and let $\disk$ be the open unit disk in $\C$.
We will denote by $f^{\circ n}$ the $n^{th}$ iterate of
the function $f$, and by $\Sphere$ a 2-sphere with two distinguished
points $0$ and $\infty$. A holomorphic exponential map will be
written as $E_\lambda(z) := \lambda \exp(z)$.

We will use the following standard concepts on exponential dynamics;
compare \cite{SZ1,SZ2}.

\begin{defn}\capt{Escaping Point}\label{def:escpt}%
For an entire holomorphic function $f$, an {\em escaping point} is a point
$z \in \C$ with $f^{\circ n}(z) \rightarrow \infty$ as $n \rightarrow \infty$; its
orbit is an {\em escaping orbit}.
\end{defn}

A holomorphic exponential map $E_\lambda$ will be called {\it
postsingularly finite} if its singular value $0$ has a finite orbit, which means that the singular
orbit is preperiodic (we use the word ``preperiodic'' in the strict sense, excluding
the periodic case; the term ``(pre)periodic'' is used to mean either ``periodic'' or ``preperiodic'').

The following discussion applies to exponential maps for which the singular orbit does not escape; only this case is of interest to us.

\begin{defn}\capt{Dynamic Ray} \label{def:dynray}%
A {\em dynamic ray of $f$} is a maximal injective curve $\gamma
\colon(0,\infty) \rightarrow {\C}$ with $\gamma(t) \to\infty$ as $t\to\infty$ so that
 $\gamma(t)$ is an escaping point for each $t \in
(0,\infty)$. The dynamic ray {\em lands at a point $a \in \C$} if
$\lim_{t\rightarrow 0} \gamma(t) = a$.
\end{defn}

In \cite[Theorem~4.2]{SZ1}, dynamic rays were defined as curves consisting of escaping points
and satisfying certain asymptotic properties (as in Theorem~\ref{thm:dynray@extadd} below).
It was shown \cite[Corollary~6.9]{SZ1} that every escaping point is either on a unique dynamic ray,
or it is the landing point of a unique dynamic ray. Every path component of the set of
escaping points is a dynamic ray \cite[Corollary~4.3]{FRS}; as such it comes
with a parametrization as an injective curve. Therefore, our definition of dynamic rays
given above coincides with the original one in \cite{SZ1} (and is easier to state).

The $E_\lambda$-image of any dynamic ray is contained in a dynamic ray; if the singular value $0$ does not escape, then the $E_\lambda$-image of every dynamic ray equals a dynamic ray.
A dynamic ray $\gamma$ is {\it periodic} if there is an $n\ge 1$
such that $\gamma((0,\infty)) \supset E_{\lambda}^{\circ
n}(\gamma(0,\infty))$ and it is {\it preperiodic} if $
E_{\lambda}^{\circ k}(\gamma(0,\infty))$ is periodic for some
$k>0$. Note that no point on a ray can be periodic or
preperiodic since it escapes, but the curve as a set can be.

\begin{defn}\capt{External Address} \label{def:extadd}%
An {\em external address} $\s$ is a sequence $\s=s_1s_2s_3\dots$ over the integers.
Let $\Sym$ be the space of all external addresses endowed with the lexicographic order, and let  $\sigma\colon \Sym\to \Sym$ be the (left) shift map.
\end{defn}

The meaning of the external address in the dynamics of the
exponential map is as follows. For an exponential map $E_\lambda$
and for each $j\in \Z$, we let
\begin{equation*}\label{eq:static}
R_j = \{z\in\C\ \colon\ -\Im(\log \lambda) - \pi + 2\pi j< \Im(z) <
-\Im(\log \lambda) + \pi + 2\pi j\},
\end{equation*}
using the convention that $-\pi < \Im(\log \lambda) \leq \pi$. On
each $R_j$, $E_\lambda$ is a conformal isomorphism onto $\C'$. The
boundaries of the strips are the set $E_{\lambda}^{-1}(\R^{-})$.
This partition of the complex plane into strips is called the {\it
static partition}. The choice of labels for the strips is so that
$0\in \ovl R_0$.

\looseness -1 A dynamic ray $\gamma$ has external address $\s$ if
for all $n\in \N$ there is $r_n \in \R^+$ so that
$E_{\lambda}^{\circ n}(\gamma(t)) \in R_{s_n}$ if $ t>r_n$ (dynamic
rays may well cross the static partition, but they do so only  for
bounded values of $t$). By construction in \cite[Theorem~4.2]{SZ1},
different dynamic rays have necessarily different external
addresses: external addresses are the analog to external angles for
polynomial rays.
Different rays may land at the same point; these rays will then of
course have different external addresses; see
Section~\ref{sec:holexpmap2extadd}. A dynamic ray is (pre)periodic
if and only if its external address is.

In \cite[Theorem~4.2 and Corollary 6.9]{SZ1}, a complete classification of escaping points and thus of dynamic rays was given. For our purposes, the following special case is sufficient:

\begin{thm}\capt{Dynamic Ray at External Address}
\label{thm:dynray@extadd}%
If the singular orbit does not escape, then for every bounded external address $\s$ there is a unique
injective continuous curve $g_{\s}\colon(0,\infty)\rightarrow \C$ (the dynamic ray at external address $\s$)
consisting of escaping points such that:
\[
\lim\limits_{t\rightarrow \infty} \Re(g_{\s}(t)) = +\infty
\]
satisfying
\[
E_{\lambda}(g_{\s}(t)) = g_{\sigma
(\s)}(F(t)),\ \forall\ t>0
\]
 and
 \begin{equation}
 g_{\s}(t) = t - \log\lambda + 2\pi i s_1 + r_{\s}(t)
 \label{Eq:DynRayAsymptotics}
 \end{equation}
with $|r_{\s}|< 2e^{-t}(|\log \lambda |+C)$ and $F(t) = e^t -1$; here $C$ is a universal constant and $\log\lambda$ denotes a branch with $|\Im(\log\lambda)|\le\pi$.
\qedd
\end{thm}

\hide{
Consider the symbol space $\Sym = \{ s_1s_2s_3\dots$: all $s_i \in
\Z\}$ with the lexicographic order and the order topology.
}
For our combinatorial classification of postsingularly finite exponential maps, we need
a few concepts from symbolic dynamics. In what
follows, terms like $t_1\t$ (with $t_1 \in \Z$ and $\t \in \Sym$) will
denote concatenation. Let $\t = t_1t_2\dots\in \Sym$, and suppose
$\t$ is not a constant sequence. Then either $\t \in (t_1\t,
(t_1+1)\t) \subset \Sym$ or $\t \in ((t_1-1)\t, t_1\t) \subset
\Sym$; denote the interval that contains $\t$ by $I_0$. For $\u\in
\Z$, define the intervals
\[
I_\u:=\{s_1s_2s_3\dots\in\Sym\colon (s_1-\u)s_2s_3\dots\in
I_0\}\,\,.
\]
\hide{$I_\u$ in the same way using $t_1 + \u$ instead of $t_1$. }
Then $\bigcup_{\u \in \Z}I_\u$ is a partition of
$\Sym\sm\{\sigma^{-1}(\t)\}$. Using this, we can define combinatorial itineraries:

\begin{defn}
\label{def:itinerary} \capt{Itinerary $\It{\s}{\t}$}%
Consider two sequences $\s,\t \in \Sym$. The itinerary of $\s$ with respect to $\t$, denoted
$\It{\s}{\t}$, is the sequence $\uu =
\u_1\u_2\u_3\dots$ over $\Z$ such that $\sigma^k(\s) \in
I_{\u_{k+1}}$ for $k \geq 0$, where the $I_\u$ are defined as above.
If $\sigma^k(\s)\in\partial I_\u$ for some $k$ (hence $\sigma^{k+1}(\s)=\t$), we leave the itinerary undefined; this case will not be needed here.
\end{defn}

In order to motivate this formal definition, consider the dynamic ray $\gamma$ at external address $\t$ and suppose it does not contain the singular value. The countably many $\El$-preimages of $\gamma$ are dynamic rays at external addresses $k\t$ for $k\in\Z$ (where again $k\t$ denotes concatenation). These preimage rays subdivide right half planes $\{z\in\C\colon \Re(z)>x\}$ (for sufficiently large $x$) into countably many components, and every unbounded component contains unbounded parts of exactly the dynamic rays $\ray$ at external addresses $\s\in I_\u$ for one particular choice of $\u\in\Z$, or equivalently those rays $\ray$ whose
itineraries $\It{\s}{\t}$ have a given first entry $\u$. The lexicographic order of addresses corresponds exactly to the vertical order of rays in their approach to $\infty$. In Section~\ref{sec:combitins}, we will show that this concept makes particular sense for postsingularly finite exponential maps.

We are now ready to state the main theorems. Together, they give a complete combinatorial coding of postsingularly finite exponential maps: we construct a map from preperiodic external addresses to postsingularly finite exponential maps. The first theorem shows that the map is well-defined and surjective, the second one measures how injective this map is and thus defines an equivalence relation on preperiodic external addresses in terms whether or not they describe the same map.

\begin{thm}\capt{Combinatorial Coding of Exponential Maps}\label{thm:main}%
For every preperiodic external address $\s$ starting with
the entry $0$, there is a unique postsingularly finite exponential
map such that the dynamic ray at external address $\s$ lands at the singular value.

Every postsingularly finite exponential map is associated in this way to a positive finite number of preperiodic external addresses starting with $0$.
\end{thm}

\begin{thm}\capt{Different Codings}\label{thm:differentcodings}%
For any two preperiodic external addresses $\s$ and
$\s'$, the following are equivalent:
\begin{enumerate}
\item
there is a postsingularly finite exponential map $\El$ so that in its dynamic plane, the dynamic rays at external addresses $\s$ and $\s'$ land at the singular value;
\label{Item:Equiv_dynrays}
\item
the parameter rays at external addresses $\s$ and $\s'$ (see Section~\ref{sec:holexpmap2extadd}) land at the same parameter $\lambda$;
\label{Item:Equiv_pararays}
\item
$\It{\s'}{\s} = \It{\s}{\s}$;
\label{Item:Equiv_itin1}
\item
$\It{\s}{\s'} = \It{\s'}{\s'}$;
\label{Item:Equiv_itin2}
\item
$\It{\s}{\s'} = \It{\s}{\s}=\It{\s'}{\s} = \It{\s'}{\s'}$;
\label{Item:Equiv_itin1&2}
\end{enumerate}
In all these cases, $\s$ and $\s'$ have equal period and equal preperiod.

If $\s$ is a preperiodic external address with preperiod $l$ and period $k$, then the itinerary $\It{\s}{\s}$ (the kneading sequence of $\s$) has also preperiod $l$ and period $k'$ dividing $k$. The exact number of external addresses which yield the same postsingularly finite exponential map is equal to $k/k'$ if $k>k'$, and it equals $1$ or $2$ if $k=k'$.
\end{thm}

The above two theorems give a complete classification of postsingularly finite exponential maps in terms of external addresses. With some more combinatorial efforts, one can turn this into a classification by {\em internal addresses} as defined in \cite{LS,RS}; in this setting, every postsingularly finite exponential map is described by a unique internal address, which is a strictly increasing sequences of positive integers for which the difference sequence is eventually periodic, and subject to a certain admissibility condition. We do not discuss this here (see the section on unicritical polynomials in \cite{BS}).

The proof of our classification result uses the main result of
Hubbard, Schleicher and Shishikura \cite{HSS}, which is an extension
of Thur\-ston's fundamental theorem on rational maps to the setting of
exponential maps: their theorem is used in the existence part of our
statement.

The global approach to our results is illustrated in the following commutative diagram. Our exponential maps are always assumed to be postsingularly finite and the external addresses to be preperiodic.

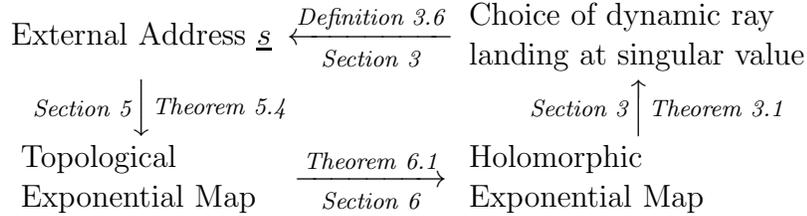
\begin{equation}
\begin{CD}
\mbox{External Address} \ \s @<{\mbox{\scriptsize\em Definition~\ref{cor:dynray2extadd}}}<{\mbox{\scriptsize\em Section~\ref{sec:holexpmap2extadd}}}<
\begin{minipage}{45mm} \setlength\baselineskip{15pt}
Choice of dynamic ray \\ landing at singular value\end{minipage}
%\mbox{Choice of dynamic ray landing at 0}
\\
@V{\mbox{\scriptsize\em Section~\ref{sec:extadd2topexpmap}}}V{\mbox{\scriptsize\em Theorem~\ref{thm:extadd2topexpmap}}}V  @A {\mbox{\scriptsize\em Section~\ref{sec:holexpmap2extadd}}} A {\mbox{\scriptsize\em Theorem~\ref{thm:prepray@singval}}} A\\
\begin{minipage}{32mm} \setlength\baselineskip{15pt}
Topological \\ Exponential Map \end{minipage}
%\mbox{Topological Exponential Map}
@>{\mbox{\scriptsize\em Theorem~\ref{thm:HSS}}}>{\mbox{\scriptsize\em Section~\ref{sec:topexp2holexp}}}>
\begin{minipage}{45mm} \setlength\baselineskip{15pt}
Holomorphic \\ Exponential Map \end{minipage}
%\mbox{Holomorphic Exponential Map}
\end{CD}\nonumber
\end{equation}
\medskip

We start our classification in Section~\ref{sec:holexpmap2extadd}:
for every postsingularly finite holomorphic exponential map, a finite
positive number of preperiodic dynamic rays lands at the singular
value; choose one such ray. Every dynamic ray has a unique
associated external address; it turns out that dynamic rays landing
at the singular value always have external addresses starting with
$0$. So far, this associates to every postsingularly finite exponential
map a preperiodic external address (this involves a choice).
In Section~\ref{sec:holexpmap2extadd}, we also discuss rays in parameter space.

In Section~\ref{sec:combitins}, we introduce some more concepts and
algorithms from symbolic dynamics which we will need in the sequel.
The main technical construction then comes in
Section~\ref{sec:extadd2topexpmap}: for every preperiodic external
address we first construct a graph with a continuous self-map and
then extend it to a branched covering of $\Sphere$ which we call a
{\em topological exponential map}. Symbolic dynamics helps us to set
things up so that there is no Thurston obstruction. Therefore, in
Section~\ref{sec:topexp2holexp} we can apply Thurston theory
(applied to exponential maps) to find an equivalent holomorphic
postsingularly finite exponential map, and again symbolic dynamics
shows that the ray at external address $\s$ lands at the singular
value. This finally shows that there is a well-defined and
surjective map from preperiodic external addresses to postsingularly
preperiodic exponential maps, so that the exponential map associated
to an address $\s$ has the property that the dynamic ray $\ray$
lands at the singular value. Finally, we investigate which external
addresses give rise to the same holomorphic exponential map, thus
describing exactly how far this map is from being injective.

\hide{
We should add that a related classification of postsingularly finite exponential maps can be given in terms of {\em internal addresses} as developed in \cite{LS}; see also the appendix in \cite{RS}: every external address gives rise in a natural way and by a simple combinatorial algorithm to an angled internal address. Then our classification simply states that two preperiodic external addresses yield the same postsingularly finite exponential map if and only if they yield the same angled internal address.
}

\section{From Exponential Map to External Address}\label{sec:holexpmap2extadd}
%Bottom right to top left arrow.++ Vlad:

In this section, we start with a postsingularly finite exponential map $\El$. We show that a preperiodic dynamic ray lands at the singular value, and associate to $\El$ the external address of the ray.

The hardest part of the work  has conveniently been done in \cite[Theorem~4.3]{SZ2}:

\begin{thm}\capt{Preperiodic Ray at Singular Value} \label{thm:prepray@singval}%
For every postsingularly finite exponential map, at least one and at most finitely many
preperiodic dynamic rays land at the singular value.
\qedd
\end{thm}

\looseness -1
There can be several dynamic rays landing at the singular value. Our classification uses the fact that all of them start with the entry $0$.

\begin{prop}\capt{External Address Starts With $0$} \label{Prop:FirstEntryZero}%
If the dynamic ray $\ray$ lands at the singular value for a postsingularly finite exponential map, then the external address $\s$ starts with $0$.
\end{prop}

This is not an obvious statement: the external address of a dynamic
ray is defined using the asymptotics for large real parts; a priori,
it seems quite possible that dynamic rays with non-zero first
entries in their external addresses could make it to the singular
value. We prove this result at the end of this section, but we will
need to introduce parameter rays (and also for other
purposes). Note that this happens in reversal of Douady's famous
principle ``you first plough in the dynamical plane and then harvest
in parameter space''. We do this the other way around (like Rempe in
\cite{Lasse}).

Similarly as dynamic rays give structure to dynamical planes,
parameter space gets a lot of structure through parameter rays; the
latter also help to understand bifurcations of exponential maps.
Just as for quadratic polynomials and the Mandelbrot set
\cite{Orsay} as well as for higher degree unicritical polynomials and Multibrot
sets \cite{Dominik}, there are deep relations between the structure
in dynamical planes and in parameter space. We follow the arguments
from \cite[Section~IV.6]{Habil}. We will need the following special
case of the main result in \cite{FS}:

\begin{prop}\capt{Parameter Rays}\label{prop:pararays}%
For every bounded sequence $\s\in\Sym$ starting with $0$, there is an injective curve $\Ray\colon(0,\infty)\to\Cstar$ in parameter space, so that for every $t>0$, the parameter $\lambda=\Ray(t)$ is the unique parameter $\lambda$ so that for $\El$, the singular value $0=\ray(t)$. These parameter rays are disjoint for different external addresses $\s$.
\qedd
\end{prop}
The general statement in \cite{FS} deals also with unbounded external addresses, but all we need here are preperiodic hence bounded addresses.

\begin{thm}\capt{Landing of Preperiodic Parameter Rays}
\label{thm:landing_para_rays}%
For every postsingularly finite exponential map $\Elz$ and every preperiodic external address, the dynamic ray $\ray$ lands at the singular value if and only if the parameter ray $\Ray$ lands at $\lambda_0$.
\end{thm}

\hide{
For every preperiodic external address $\s$, the parameter ray $\Ray$ lands at a parameter $\lambda\in\Cstar$ so that the exponential map $\El$ is postsingularly finite, and in its dynamical plane, the dynamic ray $\ray$ lands at the singular value.
Conversely, if $\El$ is a postsingularly finite exponential maps, then $\lambda$ is the landing point of all those parameter rays at preperiodic addresses $\s$ which land at the singular value in the dynamical plane of $\El$.
\end{thm}
}

\proof
Suppose $\Elz$ has the property that the dynamic ray $\ray$ at preperiodic external address $\s$ lands at the singular value; then the singular orbit for $\Elz$ is preperiodic and terminates at a necessarily repelling periodic orbit. There is then a neighborhood $U\ni\lambda_0$ in parameter space and a unique holomorphic function $z\colon U\to\C$ so that for every $\lambda\in U$, the point $z(\lambda)$ is preperiodic with $z(\lambda_0)=0$, and $z(\lambda)$ is still the landing point of the dynamic ray $\ray$. This follows from the same arguments as in \cite{GM} for the polynomial case: it suffices to know that $z(\lambda_0)$ can be extended holomorphically as a repelling preperiodic point (this is the implicit function theorem) and that for fixed potentials $t>0$, the point $\ray(t)$ depends holomorphically on $\lambda$ (and this follows from \cite[Proposition~3.4]{SZ1}). If $\lambda$ makes a small loop around $\lambda_0$, there must be at least one parameter along this loop for which $\ray$ contains the singular value $0$: during one loop of $\lambda$ around $\lambda_0$, the landing point $z(\lambda)$ must loop some number $n\neq 0$ times around $0$ (where $n$ is the local degree of the holomorphic map $\lambda\mapsto z(\lambda)$); the same is thus true for points $\ray(t)$ with very small potentials $t$. However, this is not so for large potentials $t$ because of the asymptotics in Theorem~\ref{thm:dynray@extadd}, and this proves the claim.
If $0=\ray(t)$ for $\El$, this means $\lambda=\Ray(t)$ by Proposition~\ref{prop:pararays}. Since this is true for arbitrarily small loops, $\lambda_0$ must be a limit point of $\Ray$.

Suppose that $\lambda_1\in U$ was another limit point of $\Ray$ with $z(\lambda_1)\neq 0$. For this parameter, the dynamic ray $\ray$ lands at $z(\lambda_1)$ by definition of $U$, and in particular the singular value is not on $\ray$ or on one of the finitely many rays on the forward orbit of $\ray$. Since  $z(\lambda_1)\neq 0$, and the ray $\ray$ together with its landing point form a compact set which changes continuously with $\lambda$ (again in analogy to \cite{GM}), it follows that $\lambda_1$ has a neighborhood of parameters $\lambda$ in which $0\notin\ray$. But this contradicts the assumption that $\lambda_1\in U$ was a limit point of $\Ray$. Therefore, the only limit points of $\Ray$ within $U$ can be $\lambda_0$, plus possibly finitely many further parameters $\lambda$ with $z(\lambda)=0$. The set of limit points of any ray is always connected, so $\Ray$ lands at $\lambda_0$.

Conversely, suppose $\lambda_0$ is the landing point of the parameter ray $\Ray$. Then by \cite[Theorem~3.2]{SZ2}, the dynamic ray $\ray$ lands at a repelling preperiodic point $z_0$. Similarly as above, ray and landing point are stable under perturbations. If $z_0\neq 0$, then $\lambda_0$ could not even be a limit point of $\Ray$.
\qed

The following result is stated for convenient reference.

\begin{coroll}\capt{Landing Properties of Preperiodic Parameter Rays} \label{Cor:PrepParaRays}%
Every parameter ray $\Ray$ at preperiodic external address $\s$ lands at a postsingularly finite exponential map, and every preperiodic exponential map is the landing point of a finite positive number of parameter rays at preperiodic external addresses.
\end{coroll}
\proof
This follows immediately as soon as our classification theorems are proved (we will not use it before).
\qed

\proofof{Proposition~\ref{Prop:FirstEntryZero}}
It is shown in \cite{FS} (or \cite[Corollary~3.2]{FRS}) that if the singular value escapes on a dynamic ray $\ray$, then the external address $\s$ starts with $0$ (provided dynamic rays are parametrized so that $|\Im\log\lambda|<\pi$). If a parameter ray lands at a postsingularly finite exponential map, then rays and their parametrization change continuously.
\qed

It might seem that the statement of Proposition~\ref{Prop:FirstEntryZero} makes sense only once a branch of $\log\lambda$ is chosen, which is not a dynamically well-defined quantity. However, this is not so: the proposition says that independently of any choice of branch, and any choice of labels of strips defining external addresses, any dynamic ray landing at $0$ has asymptotic imaginary part in $(-\pi,\pi)$.

\begin{defn}\label{cor:dynray2extadd}\capt{External Address of $E_\lambda$}%
Let $E_\lambda$ be a postsingularly finite holomorphic exponential
map. Then we associate to $\El$ the external address of a dynamic
ray $\ray$ which lands at $0$ (this may involve a choice).
\end{defn}

\section{Symbolic Dynamics and Kneading Sequences}
\label{sec:combitins}

In Section \ref{sec:extadd2topexpmap} we aim
to construct a topological exponential map $f$ in which we encode
all the combinatorial information of a given external address
$\s$. In order to do this, we need a few more concepts from symbolic dynamics.

In Definition~\ref{def:itinerary}, we defined the space $\Sym$ of external addresses and, for
every pair of sequences $\s,\t\in\Sym$, the itinerary $\It{\s}{\t}$ of $\s$ with respect to $\t$.
Of special importance is the itinerary of a sequence with itself: the kneading sequence.

\begin{defn}
\label{def:kneadingseq} \capt{Kneading Sequence}%
For a sequence $\s \in \Sym$ we call $K(\s):=\It{\s}{\s}$ the
{\em kneading sequence} of $\s$.
\end{defn}

The methods of symbolic dynamics and the concept of itineraries are especially useful
for those exponential maps $\El$ for which a dynamic ray $\ray$ lands at the singular value:
in particular, if the singular orbit is preperiodic (the main case of interest to us), then
by Theorem~\ref{thm:prepray@singval} there are one or several dynamic rays
at preperiodic external addresses landing at $0$
(see also \cite[Section~4]{SZ2} for a discussion of several other cases with similar properties).
In this case, the countably many $\El$-preimages of $\ray$ partition all of $\C$ and
form what we call a {\em dynamic partition}. The components in this partition are
translates of each other by $2\pi i\Z$; the imaginary parts of
any component are in general unbounded (usually, the ray $\ray$
spirals into its landing point $0$).

There is always a unique component, called $I_0$, which contains the singular value, and its vertical translate by $2\pi j$ is called $I_j$ for $j\in\Z$. If $z\in\C$ is a point whose orbit is disjoint from $\ray$, then we define the {\em itinerary} of $z$ (with respect to the ray $\ray$) as the sequence of component labels visited by the orbit of $z$.

We call this new partition the {\em dynamic partition} (as opposed to the static partition introduced in Section~\ref{sec:defns}). The dynamic partition has the advantage that each dynamic ray is
completely contained in one component (unless it is one of the rays forming the partition boundary),
and all points on the ray and its possible landing point have the same itinerary. In fact,
the itinerary of all points on the ray at external address $\t$ is $\It{\t}{\s}$, and the itinerary of the singular value (or of any point on any ray landing at the singular value) is the kneading sequence $K(\s)$. The following result is shown in \cite[Proposition~4.4]{SZ2}.

\begin{lemma}\capt{Itinerary of Landing Points and Rays}
\label{lemma:itin&pointland}%
For postsingularly finite exponential maps, no two (pre)periodic
points have the same itinerary, and a (pre)periodic dynamic ray
lands at a given periodic or preperiodic point if and only if ray
and point have the same itinerary. In particular, two (pre)periodic
dynamic rays land together if and only if they have the same
itineraries.
\qedd
\end{lemma}

Note that for the dynamic partition, unlike the static partition, several rays may have the same itinerary.

The following simple algorithm illustrates the close relation
between external addresses and itineraries and shows in particular
that their entries differ at most by $1$
(up to simultaneous translation of all entries by the same integer).

\begin{algo}
\label{algo:itinerary} \capt{Construction of $\It{\s}{\t}$}
Given external addresses $\s, \t \in\Sym$ so that $\t$ is non-constant and $\sigma^n(\s)\neq\t$ for all $n\ge 0$. Then $\uu := \It{\s}{\t}$ can be constructed as follows.
\begin{enumerate}
\item
For $n \geq 1$ define  $\delta_n := \left\{
\begin{array}{ll} -1 & \mbox{ if } \sigma(\t) > \t \mbox{
and } \sigma^{\circ
n}(\s) < \t\\
\ \ 1 & \mbox{ if } \sigma(\t) < \t \mbox{ and }
\sigma^{\circ
n}(\s) > \t\\
\ \ 0 & \mbox{ otherwise.} \end{array} \right.$
\item
Construct $\uu = \u_1 \u_2\dots$ as \; $\u_n = s_n - t_1 + \delta_n$ \,\, .
\end{enumerate}
\end{algo}

\proof
It suffices to show that the first entry $\u_1$ in $\It{\s}{\t}$ is correct: the $n$-th entry equals by definition the first entry in $\It{\sigma^{n-1}(\s)}{\t}$. Adding an integer $k$ to the first entry $s_1$ of $\s$ will add $k$ to $\u_1$, so we may assume that $s_1=t_1$.

Suppose first that $\s > \t$. Then $\u_1=0$ unless there is a
preimage of $\t$ in $(\t,\s)$; but since $s_1=t_1$, this is
equivalent to the condition $\t\in(\sigma(\t),\sigma(\s))$ or
$\sigma(t)<t<\sigma(\s)$; and exactly in this case, $\u_1=1$.
Similarly, if $\s<\t$, then $\u_1=0$ unless
$\sigma(t)>t>\sigma(\s)$, and exactly in that case, $\u_1=-1$. \qed

Notice that we have $\sigma^n(\s) = \s$ for some $n$ if and only if
$\s$ is periodic. Algorithm \ref{algo:itinerary} therefore works for
computing the kneading sequences of the preperiodic addresses that
we are interested in. It will prove to be useful when we are trying
to recover the external address from our constructed exponential
map.

\section{The Topological Exponential Map}
\label{sec:extadd2topexpmap}

In this section, we will start with a combinatorial object (external
address) and turn this into a topological object (a postsingularly
preperiodic topological exponential map). In the next section, we
make the step from topology to the complex structure and find,
whenever possible, a holomorphic exponential map which is equivalent,
in a sense defined by Thurston,
to the given topological exponential map.

{\sc Convention.} All homeomorphisms and coverings in this paper
will be orientation preserving.

\begin{defn}
\capt{Topological Exponential Map} \label{def:topexpmap}%
A universal cover $f\colon (\Sphere\sm \{\infty\}) \to (\Sphere\sm
\{\infty, 0\})$ is called a topological exponential map. It is
called postsingularly finite if the orbit of 0 is finite, hence
preperiodic. The postsingular set is $P_f := \bigcup_{n \geq 0}
f^{\circ n} (0) \cup \{\infty\}$.
\end{defn}

If a topological exponential map is holomorphic, then it is conformally conjugate
to an exponential map  $z\mapsto\El$.

\begin{defn}
\capt{Thurston Equivalence} \label{def:thurstequiv}%
Two postsingularly finite exponential maps $f$ and $g$ with
postsingular sets $P_f$ and $P_g$ are called Thurston equivalent
if there are two homeomorphisms $\phi_1, \phi_2 \colon\Sphere
\to \Sphere$ with $\phi_1|_{P_f} = \phi_2|_{P_f},\ P_g =
\phi_1(P_f) = \phi_2(P_f)$ and $\phi_1(\infty) = \phi_2(\infty) =
\infty$ such that
\[
\phi_1 \circ f = g \circ \phi_2 \quad \ \mbox{on}\
\Sphere\sm\{\infty\}
\]
and $\phi_1$ is homotopic (or equivalently isotopic) to $\phi_2$ on $\Sphere$ relative to $P_f$.
\end{defn}

Our goal will be to find, for every postsingularly finite
topological exponential map, a postsingularly finite holomorphic
exponential map which is Thurston equivalent. This is not always
possible. In the case of rational mappings, Thurston \cite{DH}
determined that this is impossible if and only if there is what is now called
a Thurston obstruction; see also \cite{BFH,HS}. The extension
of this result to the case of exponential maps was done in
\cite{HSS}: in this case the possible obstructions have a much
simpler form, called degenerate Levy cycles.

\begin{defn} \capt{Essential Curves and Levy Cycle} \label{def:levy}%
Let $f$ be a topological exponential map with postsingular set $P_f$.
A simple closed curve $\delta\subset\Sphere\sm P_f$ is called {\em essential} if both components of $\Sphere\sm\delta$ contain at least two points in $P_f$.

Suppose there exist disjoint essential simple closed curves $\delta_0,\dots, , \delta_k =
\delta_0 $ such that for each $i = 0, \dots, k-1$, $\delta_i$ is
homotopic relative $P_f$ to one component $\delta'$ of
$f^{-1}(\delta_{i+1})$ and $f\colon \delta' \rightarrow \delta_{i+1}$
has degree 1. Then $\Lambda=\{\delta_0,\delta_1,\dots, , \delta_k =
\delta_0  \}$ is called a {\em Levy cycle}.
\end{defn}

Essential curves are important for the following reason: a simple closed curve $\delta\in\Sphere\sm
P_f$ is essential if and only if, for every homeomorphism
$\phi\colon\Sphere\to\Cbar$, there is a lower bound of lengths (with
respect to the hyperbolic metric of $\Cbar\sm \phi(P_f)$) of simple
closed curves homotopic to $\phi(\delta)$ relative to $\phi(P_f)$.

This section will be concerned with proving the following theorem:

\begin{thm}\capt{External Address Yields Topological Exponential Map}
\label{thm:extadd2topexpmap}%
Let $\s$ be a preperiodic external address. Then there exists a postsingularly finite topological exponential map $f$ with the following properties:
\begin{itemize}
\item $f$ has a preperiodic injective curve $\gamma\colon(0, \infty)\to\Sphere$
connecting $0$ to $\infty$;
\item
$\gamma$ has  itinerary $\It{\s}{\s}$ with respect to the partition defined by $f^{-1}(\gamma)$,
\item the vertical order of the rays $f^{\circ n}(\gamma)$ coincides with the lexicographic order of the shifts of $\s$,
\item $f$ does not admit a Levy cycle.
\end{itemize}
Any two such postsingularly finite topological exponential maps for the same external address $\s$ are Thurston equivalent to each other.
\end{thm}
Note that any injective curve connecting the singular value $0$ to $\infty$ has countably many disjoint preimages under any topological exponential map, and this allows us to define a dynamic partition and thus dynamical itineraries just like for holomorphic exponential maps for which a dynamic ray lands at $0$.

As always, the curve $\gamma$ should be preperiodic as a set; its points need not be (except the endpoint). The preimage $f^{-1}(\gamma)$ is disjoint from all rays $f^{\circ n}(\gamma)$. Let $\gamma'$ be the unique component of $\gamma\cap\disk$ starting at $0$ and let $p'$ be any component of $f^{-1}(\gamma')$. Then the rays $f^{\circ n}(\gamma)$, as well as $p'$, are disjoint curves to $\infty$ and have a well-defined cyclic order. Removing $p'$ induces a linear order among all rays, and this is the vertical order specified by the theorem; it does not depend on the choice of $p'$.

\subsection{The Graph Map}
\label{Sub:GraphMap}

Similarly as for polynomials in \cite{BFH}, we start by constructing
an undirected graph $\Gamma$ that encodes the combinatorial
information given by $\s$. An important difference is that our
graph is infinite.

We will construct an infinite topological graph $\Gamma$ and later
embed it into $\Sphere$. Start with two vertices
$\Gamma=\{e_\infty,e_{-\infty}\}$. For each $n \in \Z$, add disjoint
edges $p_n$ joining $e_\infty$ to $e_{-\infty}$. Let $k$ and $l$ be
the length of preperiod and period of $\s$  respectively. Add
vertices $e_1, \dots, e_{k+l}$ to $\Gamma$, and for each $e_n$, add
an edge $\gamma_n$ connecting $e_n$ and $e_\infty$, so that all
edges are disjoint and all vertices are disjoint from each other and
from all edges.

We will embed $\Gamma\sm\{e_{-\infty}\}$ into $\Sphere$ and define a
graph map $\tilde{f}$ from the  embedded graph to itself. By
extending $\tilde{f}$ to a map $f\colon\Sphere\sm\{\infty\}
\rightarrow\Sphere\sm\{0,\infty\}$, we will obtain the desired
topological exponential map.

It is straightforward to embed $\Gamma\sm\{e_{-\infty}\}$ (as it has
been constructed up to here) into $\Sphere$ in a reasonable way,
define a graph map, and extend it to a topological exponential map
$f$ on all of $\Sphere$, so that $f$ satisfies the first three
properties of Theorem \ref{thm:extadd2topexpmap}. The hard part is
to make sure that $f$ will not admit a Levy cycle. The following
lemma tells us when to expect a Levy cycle:

\begin{lemma}\capt{Levy Cycle and Itineraries}
\label{lemma:levy&itin}%
Consider a topological exponential map that satisfies the first two
properties in Theorem \ref{thm:extadd2topexpmap}. Then two
or more postsingular points are surrounded by the same curve in a
Levy cycle if and only if they are all periodic, and they have the
same itinerary with respect to the partition consisting of preimages
of the ray landing at the singular value.
\end{lemma}

\proof
Suppose that two or more postsingular points are surrounded by a simple closed curve $\delta$ in a degenerate Levy cycle. Note first that $\delta$ cannot surround the singular value $0$: otherwise, the preimage of $\delta$ would not contain any simple closed curve. After homotopy, we may thus assume that $\delta$ does not intersect the ray $\gamma_1$ connecting $0$ to $\infty$.

Taking preimages, no preimage curve $\delta'$ of $\delta$ can intersect the partition boundary, hence all postsingular points surrounded by $\delta'$ have the same first entries in their itineraries. Note that the number of postsingular points surrounded by $\delta'$ cannot be greater than that for $\delta$; this number could be smaller, depending on which branches of preimages are chosen. However, since $\delta$ is part of a degenerate Levy cycle and hence periodic (up to homotopy), the number of surrounded postsingular points must remain constant. It follows that all postsingular points surrounded by $\delta$ have the same periodic itinerary and are hence periodic points.

To prove the converse, assume that the periodic postsingular points $e_{i_1},\
e_{i_2},\ \ldots,\ e_{i_k}$ have the same itinerary with respect
to the dynamic partition of the plane. Surround these points (but no other postsingular points) by a simple closed curve $\delta$; this curve is automatically essential. Note that there are in general infinitely many homotopy classes of curves $\delta$ relative to the postsingular points. However, there is only a single homotopy class if we require that $\delta$ must not intersect the ray into any point which is not surrounded: the complement in $\Sphere\sm\{\infty\}$ of all these rays is simply connected.
\hide{, and any simple closed curve surrounding all these points is contractible in $\C$ to the point $\infty$.}

Since $\delta$ does not intersect the ray $\gamma_1$, and all surrounded points are periodic and have the same itinerary, there is one preimage component of $\delta$ which surrounds all periodic preimages of the surrounded points, and again it does not intersect the rays into those points that it does not surround. Repeating this argument for one period of the itinerary, we obtain another curve which surrounds the same points as $\delta$ in the complement of the remaining rays, so this curve is homotopic to $\delta$ and we have a degenerate Levy cycle.
\qed

\begin{coroll}\capt{Levy Cycle and Itineraries}\label{coroll:levy&itin}%
A topological exponential map as in Lemma \ref{lemma:levy&itin}
admits a Levy cycle containing a curve surrounding the points $e_i,
e_j$ if and only if the kneading sequence $\uu=\It{\s}{\s}$ of the
curve $\gamma$ landing at the singular value has the property that
$\u_i\u_{i+1}\u_{i+2} {\dots} = \u_j\u_{j+1}\u_{j+2} \dots$.
\end{coroll}

\proof
If $\uu=\u_1\u_2\dots$ is the kneading sequence of the curve $\gamma$ landing at the singular value, then the itineraries of the points $e_i=f^{\circ(i-1)}(0)$ and $e_j$ are the appropriate shifts of the kneading sequence. Equality of itineraries can hold only if both points have equal periods and preperiods; since they are on the same preperiodic orbit, this implies that they can have identical itineraries only if they are periodic.
\qed

The Levy cycle obstruction warns us that rays $\gamma_i$ and
$\gamma_j$ should really land together at a common point $e_i=e_j$. In
order to solve this problem, we will simply ``glue'' points $e_i$
and $e_j$ in the graph: Define the equivalence relation
\begin{equation}
e_i \sim e_j :\Longleftrightarrow \u_i\u_{i+1}\u_{i+2}\dots = \u_j\u_{j+1}\u_{j+2} \dots
\nonumber
\end{equation}
on points $e_k \in \Gamma$. Redefine $\Gamma$ as $\Gamma/\sim$

After embedding into $\Sphere$, the new quotient graph $\Gamma$ will have the property that no two different vertices have the same itinerary, so there can be no Levy cycle.

In order to check that the graph can be embedded, we have to verify the following unlinking property: it never happens that there are four external addresses $\s_1<\s'_1<\s_2<\s'_2$ so that $\s_1$ and $\s_2$ have the same itinerary $\uu$, and also $\s'_1$ and $\s'_2$ have the same itinerary $\uu'\neq\uu$.
Suppose by contradiction that this problem does occur. Then
(possibly after replacing all four addresses with the same shift), we may assume that $\u_1=\u'_1$ but $\u_2\neq\u'_2$, where $\uu=\u_1\u_2\u_3\dots$ and $\uu'=\u'_1\u'_2\u'_3\dots$. Without loss of generality, assume that $\u_2\neq 0$ (the argument is symmetric in $\uu$ and $\uu'$). Then $\sigma^{-1}(I_{\u_2})\cap I_{\u_1}$ consists of a single interval which must contain $\s_1$ and $\s_2$, but it cannot contain $\s'_1$ or $\s'_2$. This proves the unlinking property.

After preparing our graph so that no Levy cycle can emerge, we embed $\Gamma$ first into
$\Cbar$ and then into $\Sphere$. The  complex structure on $\Cbar$ does not play any role, but it allows us to describe the construction more easily.

Define the embedding $\psi\colon \Gamma \rightarrow \Cbar$ as
an injective continuous map as follows:
\begin{enumerate}
  \item \label{i}
  First let $\psi(e_1) = 0$ and $\psi(e_\infty)= \infty$.
  \item \label{ii}
  Since the $\psi$-image of $\gamma_1$ must be a curve from $0$ to $\infty$, let $\psi(\gamma_1)=\R^+$.
  \item \label{iii}
  For every $n \in \Z$ let $\psi(p_n)$ be the straight
horizontal line with imaginary part $(2n - 1)\pi$.
\item\label{iv}
The images of the edges $p_n$ of $\Gamma$ define a partition of $\C$ into
strips $\Delta_n$. We label the strips in vertical order so that $\Delta_n$
denotes the strip containing $2n\pi i$ (so $0 \in \Delta_0$).
\item\label{v}
Now we are ready to bring in the combinatorial information stemming
from the kneading sequence $\uu$. For each $n=1,2,\dots,l+k$, let
$\psi(e_n)$ be a point in strip $\Delta_{\u_n}$ (taking into account
that certain points $e_n$ might be identified) and let
$\psi(\gamma_n)$ be a curve in $\Delta_{\u_n}$ connecting
$\psi(e_n)$ to $+\infty$. For simplicity, assume that imaginary
parts of $\psi(\gamma_n)$ are eventually constant. Choose these
curves so that they are disjoint from each other and from all
endpoints except their own, and choose the eventually constant
imaginary parts so that $\im(\psi(e_n)) < \im(\psi(e_m))$ if and
only if $\sigma^{n-1}(\s) < \sigma^{m-1}(\s)$ for each $n\neq m$;
this ensures that the rays respect the order prescribed by the
external address.

We have to justify that this can be done consistently if several endpoints are identified: it can never happen that two curves $\gamma_{n}$ and $\gamma_{n}$ have a common endpoint and separate $\C$ into two complementary components which both contain a curve $\gamma_{n'}$ and $\gamma_{m'}$ that should have a common endpoint different from $e_{n}=e_{m}$. This is exactly the unlinking property.
\end{enumerate}

We now map $\Cbar$ homeomorphically to $\Sphere$, mapping $0,\infty\in\Cbar$ to the two corresponding marked points $0,\infty\in\Sphere$. From now on, we view $\psi$ as a map from $\Gamma$ to $\Sphere$ and denote $\Gamma_{\Sphere}:=
\psi(\Gamma) \subset\Sphere$. For simplicity, we write $e_n$ for $\psi(e_n)$,
$\gamma_n$ for $\psi(\gamma_n)$ and $p_n$ for $\psi(p_n)$ (the vertices, the rays, and the partition boundaries). Since from now on, we only work with the embedded graphs, no confusion can arise.

\remark
In the construction of $\Gamma_{\Sphere}$, the only combinatorial information coming from $\s$ are the kneading sequence $\itin{u}$ and the lexicographic (vertical) order of the set of external addresses: $\{\s,\sigma(\s),\ldots,\sigma^{\circ (l+k-1)}(\s)\}$. We will see in Algorithm~\ref{algo:kneading&order->extaddr} that this information gives us back a unique external address $\s$, normalized so that its first entry is $0$.

We can now define a graph map $\tilde{f}\colon\Gamma_{\Sphere}\to \Gamma_{\Sphere}$ such that $\tilde{f}(\infty) = \infty$, as well as $\tilde{f}(e_j) = e_{j+1}$ and $\tilde{f}(\gamma_j) = \gamma_{j+1}$ for all $j$ (counting indices modulo the period, so that $\tilde f(\gamma_{l+k})=\gamma_{l+1}$). Furthermore, for all $k \in \Z$, define $\tilde{f}(p_k) = \gamma_1 $.

Observe that $\tilde{f}$ is neither surjective nor injective, and
that any two glued points have glued images, so the graph map
respects the gluing. Under the map $\tilde f$, the orbit of $0$ is
necessarily preperiodic. The graph map $\tilde f$ is continuous
everywhere except at $\infty$; this will not affect the extended map
$f$, which is defined on $\Sphere\sm\{\infty\}$.

\subsection{Extension of the Graph Map}
\label{Sub:ExtensionGraphMap}

In order to prove Theorem~\ref{thm:extadd2topexpmap}, we need to do three things: we need to show that the graph map can be extended to a topological exponential map (Lemma~\ref{lemma:extgraphmap}), we need to show that it satisfies the conditions given in  the theorem (Lemma~\ref{lem:embedding_satisfies}), and we need to prove the uniqueness claims (Proposition~\ref{prop:thurston-uniqueness}). The second part is the easiest, so we do it first.

\begin{lemma}\capt{Every Embedding Satisfies Conditions on Rays}
\label{lem:embedding_satisfies}%
If $f\colon\Sphere\sm\{\infty\}\to\Sphere\sm\{0,\infty\}$ is a topological exponential map which extends $\tilde{f}\colon\Gamma_{\Sphere}\to \Gamma_{\Sphere}$, then $f$ satisfies the following properties:
\begin{itemize}
\item
the curve $\gamma_1$ connects $0$ to $\infty$ and is preperiodic under $f$;
\item
with respect to the dynamic partition induced by preimages of $\gamma_1$, the curve $\gamma_1$ has itinerary $\s$;
\item
the vertical order of the rays $\gamma_n=f^{\circ(n-1)}(\gamma_1)$ coincides with the lexicographic order of the $\sigma^{n-1}(\s)$;
\item
$f$ does not admit a Levy cycle.
\end{itemize}
\end{lemma}
\proof
The first three properties follow directly from the construction. The last property follows from Corollary~\ref{coroll:levy&itin} because $f$ can have a Levy cycle only if there are two different points $e_i,e_j$ with identical itineraries, but such points have been glued together.
\qed

Our next project is the extension of the graph map to a topological exponential map. In order to ensure that the graph map can be extended to a neighborhood of every vertex, we need the following lemma.

\begin{lemma}\capt{$\tilde{f}$ Preserves Cyclic Order}
\label{lemma:cyclicorder}%
The cyclic order at $\infty$ of three or more rays with the same
itinerary (landing at the same point) is preserved by the graph map
$\tilde f$.
\end{lemma}
\proof Suppose the endpoints $e_n$, $e_{n'}$ and $e_{n''}$ are
identified, so that the three rays $\gamma_{n},
\gamma_{n'},\gamma_{n''}$ land at a common point. The vertical order
of these rays coincides with the lexicographic order of the
addresses $\sigma^{n-1}(\s)$, $\sigma^{n'-1}(\s)$ and
$\sigma^{n''-1}(\s)$; suppose without loss of generality that
$\sigma^{n-1}(\s)<\sigma^{n'-1}(\s)<\sigma^{n''-1}(\s)$. Since the
first entries in the itineraries of these three external addresses
coincide, there is a $k\in\Z$ so that
\[
k\s < \sigma^{n-1}(\s)<\sigma^{n'-1}(\s)<\sigma^{n''-1}(\s) < (k+1)\s
\]
(where $k\s$ and $(k+1)\s$ denote adjacent preimages of $\s$ under
the shift). But on every interval $(k\s,(k+1)\s)$, the shift map is
injective and preserves the cyclic order. \qed

One key construction is Alexander's trick (compare e.g.\ \cite{BFH}):
{if $f\colon\S^1 \rightarrow \S^1$ is an orientation-preserving
homeomorphism, then there exists an orientation-preserving
homeomorphism $\hat{f}\colon\overline{\D} \rightarrow \overline{\D}$
such that $\hat{f}|_{\partial \overline{\D}} = \hat{f}|_{\S^1} = f$}:
one such extension is given by
$\hat{f}(r e^{i \theta}) := r f(e^{i \theta})$.

\hide{
One key construction is the following; for a proof, see for example \cite{BFH}.
\begin{lemma}\label{lemma:alexander}\capt{Alexander's Trick}
If $f\colon\S^1 \rightarrow \S^1$ is an orientation-preserving
homeomorphism, then there exists an orientation-preserving
homeomorphism $\hat{f}\colon\overline{\D} \rightarrow \overline{\D}$
such that $\hat{f}|_{\partial \overline{\D}} = \hat{f}|_{\S^1} = f$.
Moreover, $\hat{f}$ is unique up to isotopy rel $S^1$.
\qedd
\end{lemma}
\noindent
In fact,  $\hat{f}(r e^{i \theta}) := r f(e^{i \theta})$ provides such an extension $\hat f$.
}

Moreover if $f,g\colon\S^1 \rightarrow \S^1$ are isotopic rel some
finite number of points in $\S^1$, then by extending the isotopy
to the entire disk one gets an isotopy (rel the same points) between the extensions $\hat{f}$ and
$\hat{g}\colon\overline{\D} \rightarrow \overline{\D}$.

\begin{lemma}\capt{Extension of Graph Map}\label{lemma:extgraphmap}%
The graph map $\tilde{f} \colon \Gamma_{\Sphere} \rightarrow
\Gamma_{\Sphere}$ can be extended to a topological exponential map $f \colon
\Sphere\setminus\{\infty\}\rightarrow\Sphere \setminus \{0,\infty\}$.
\end{lemma}

\proof For convenience, let us adopt the convention that the edges
$p_n$ and $\gamma_m$ contain their endpoints $\infty$ and $e_m$.
Notice that for $n\in\Z$, the set $\Sphere\sm(p_n\cup p_{n+1})$ is
the union of three disjoint open topological disks; among them,
$\Delta_n$ is the disk which does not intersect any $p_k$ for
$k\notin\{n,n+1\}$.

We shall first construct continuous maps $f_n\colon(\ovl \Delta_n\sm\{\infty\})
\to(\Sphere\sm\{\infty\})$ for each $n$, such that $f_n (p_n) = f_n(p_{n+1}) = \gamma_1$ as above, and $f_n(\Gamma_{\Sphere}\cap \Delta_n) = \tilde{f}(\Gamma_{\Sphere}\cap
\Delta_n)$. Moreover, the restriction of $f_n$ to $\Delta_n$ is an orientation preserving homeomorphism onto $\Sphere\sm\gamma_1$ which extends continuously to the boundary.

We distinguish four possible cases  for $\Delta_n$:
\begin{enumerate}
\renewcommand{\labelenumi}{(\roman{enumi})}
\item
$\Delta_n\cap\Gamma_{\Sphere}=\emptyset$;
\item
$\Delta_n\cap\Gamma_{\Sphere}=\gamma_m$ for some $m\in\{1,2,\dots,l+k\}$;
\item
$\Delta_n\cap\Gamma_{\Sphere} = \gamma_{m_1} \cup\ldots \cup \gamma_{m_i}$, for some $m_1,\dots,m_i\in\{1,2,\dots,l+k\}$ such that all endpoints $e_{m_1},\ldots e_{m_i}$ are glued;
\item
The general case: there may be combinations of several instances of case (ii) and (iii) on one domain $\Delta_n$.
\end{enumerate}

Case (i) is almost literally Alexander's trick: the graph map $\tilde f$ prescribes the boundary values on the topological disk $\Delta_n$ for the map $f_n$. The only problem is that the point $\infty$ appears twice on the boundary, mapping to $0$ on the left and to $\infty$ on the right. This causes no problem  (we do not define $f_n$ or $f$ on $\infty$).

In Case (ii) the idea of the construction is the same; we need to find a homeomorphism $f\colon\Delta_n\sm\gamma_m\to\Sphere\sm\gamma_{m+1}$ which coincides with the prescribed boundary values on $p_n$, $p_{n+1}$, $\gamma_m$ and $\infty$ (at two sides). Note that $\Delta_n\sm\gamma_m \approx \Sphere\sm\gamma_{m+1}\approx\disk $. In this case, the curve $\gamma_m$ along with its endpoints occurs twice on the boundary, and $f_n$ is already defined on $\gamma_m$ via $\tilde f$.

Now we treat Case (iii). The set $\Delta_n\sm(\gamma_{m_1} \cup \ldots \cup \gamma_{m_i})$ consists of $i$ domains, each of which is homeomorphic to $\disk$; the same is true for $\Sphere \setminus (\gamma_1 \cup \gamma_{m_1 +1} \cup \ldots\gamma_{m_i +1})$ (note that the endpoint $e_1=0$ is never glued with any other endpoint because its itinerary has longer preperiod than all others). As in Case (i), we want to extend the map $f$ along its prescribed boundary values to the appropriate image domains. In order for this to be possible, we need to assure that for each of the $i$ topological disks, the boundary rays map to the boundary rays of an appropriate image domain. Let $a$ be the landing point of the rays $\gamma_{m_1}, \ldots, \gamma_{m_i}$; then we need to make sure that the cyclic order of these rays at $a$ coincides with the cyclic order of the image rays at $f(a)$. Our construction lets rays land together iff they have identical itineraries, so this fact is assured by Lemma~\ref{lemma:cyclicorder}. Therefore, the map $f_n$ can be defined on $\Delta_n$ as well, and it is again a homeomorphism $f_n\colon\Delta_n\to\Sphere\sm\gamma_1$.

Finally, the general case (iv) can incorporate several rays, and several groups of rays, within the same domain $\Delta_n$. First observe that the vertical order of all rays $\gamma_{m_i}\subset\Delta_n$ and  $p_n,p_{n+1}\subset\partial\Delta_n$ is compatible with the cyclic order of the image rays near infinity: the vertical order of the rays is determined by their external addresses, and all external addresses within one domain $\Delta_n$ belong to one interval on which the shift map is injective. Removing from $\Delta_n$ all rays, we obtain finitely many connected components which can be treated separately because they have compatible boundary values. For every connected component, the claim follows by combining the ideas from the previous cases.

Having constructed the maps $f_n$ on all strips, we define a continuous map $f\colon\Sphere\sm\{\infty\}  \rightarrow\Sphere\sm\{0,\infty\} $, by $f|_{\Delta_n} := f_n$. The values on the boundaries of the $\Delta_n$ (except at $\infty$) match since here they coincide with the graph map $\tilde f$. We obtain a universal cover $f\colon\Sphere \setminus \{\infty\} \to \Sphere\setminus\{0,\infty\} $.
\qed

\subsection{Thurston-Uniqueness}
\label{Sub:ThurstonUnique}

Of course, the construction of the extension is not unique; however,
we have a uniqueness result in the following sense:

\begin{prop}\capt{Thurston-Uniqueness of Topological Exponential Map}
\label{prop:thurston-uniqueness}%
Suppose that $f_1,f_2\colon(\Sphere\sm\{\infty\})\to(\Sphere\sm\{0,\infty\})$ are topological exponential maps satisfying the four itemized properties of Theorem~\ref{thm:extadd2topexpmap} for the preperiodic external address $\s$.
Then $f_1$ and $f_2$ are postsingularly finite and Thurston equivalent.
\end{prop}

We start with the following lemma.

\begin{lemma}\capt{Uniqueness of Graph Map up to Homotopy}
\label{lem:uniqueness_graph_map}%
Suppose that $f\colon(\Sphere\sm\{\infty\})\to(\Sphere\sm\{0,\infty\})$ is a postsingularly finite topological exponential map satisfying the four itemized properties of Theorem~\ref{thm:extadd2topexpmap} for the preperiodic external address $\s$. Then the restriction of $f$ to the collection of curves $\bigcup_{n\ge -1}f^{\circ n}(\gamma)$ yields a graph map $\Gamma_{\Sphere}$ which is homotopic relative to $P_f$ to any graph map as constructed above.
\end{lemma}
With $\bigcup_{n\ge -1}f^{\circ n}(\gamma)$, we mean the countably many preimage curves $f^{-1}(\gamma)$, together with the finitely many curves on the forward orbit of $\gamma$.

\proof
First we find a homotopy of $\Sphere$ relative $P_f$ sending the preimages $f^{-1}(\gamma)$ to the edges $p_n$ of $\Gamma_{\Sphere}$; this is possible because the points in $P_f$ have prescribed itineraries. Then we can find homotopies within each complementary domain $\Delta_n$ of $\Sphere\sm \bigcup_n p_n$ to match the graphs within each $\Delta_n$; this is possible because
\begin{itemize}
\item
the itinerary prescribes which legs $\gamma_m$ are within which $\Delta_n$,
\item
the order at $\infty$ of the different legs within the same $\Delta_n$ is in both cases prescribed by the lexicographic order of appropriate shifts of $\s$, and finally
\item
the same endpoints of legs are identified because of the non-existence of Levy cycles: Corollary~\ref{coroll:levy&itin} says which endpoints are identified in the topological exponential map, and the same endpoints are identified by construction in our graph map.
\qed
\end{itemize}

\begin{lemma}\capt{Uniqueness of Graph Map Extension up to Homotopy}\label{lem:uniqueext}%
Let $\tilde{f}\colon \Gamma_{\Sphere} \rightarrow \Gamma_{\Sphere}$ be a
graph map and $f_1,f_2\colon\Sphere \setminus \{\infty\} \rightarrow
\Sphere\setminus\{0,\infty\}$ two different graph map extensions as
in Lemma \ref{lemma:extgraphmap}. Then the maps $f_1$ and $f_2$
are Thurston equivalent.
\end{lemma}

\proof
As in \cite{BFH}, we will define a map $\varphi\colon\Sphere \rightarrow
\Sphere$ such that $f_1 = f_2 \circ \varphi$ and $\varphi$ is a
homeomorphism homotopic to the identity rel the vertices of
$\Gamma_{\Sphere}$. The construction of $\varphi$ is easy. Each
component $U$ of $\Sphere \setminus \Gamma_{\Sphere}$ is homeomorphic to
a disk, and by construction $f_1$ and $f_2$ coincide on their boundaries. Since both $f_i$ restricted to these disks are homeomorphisms onto their images, they are homotopic to each other relative to the boundary; compare the remark after Lemma~\ref{lemma:alexander}. Since the postsingular set is contained in $\Gamma_{\Sphere}$, $f_1$ and $f_2$ are Thurston equivalent in the sense of Definition~\ref{def:thurstequiv}.
\qed

\proofof{Proposition~\ref{prop:thurston-uniqueness}}
Since the curve $\gamma$ is preperiodic as a set, and one of its endpoints is the singular value $0$, it follows that $f_1$ and $f_2$ are postsingularly finite. By Lemma~\ref{lem:uniqueness_graph_map}, the graph maps of $f_1$ and $f_2$ coincide up to homotopy relative to the postsingular set; the homotopy of graph maps extends to a homotopy of $\Sphere$ because the graph maps are embedded into $\Sphere$ in the same way, and by Lemma~\ref{lem:uniqueext}, we obtain a Thurston equivalence between $f_1$ and $f_2$.
\qed

\proofof{Theorem~\ref{thm:extadd2topexpmap}}
We have first constructed a graph $\Gamma_{\Sphere}\subset\Sphere$ and a graph map $\tilde f\colon\Gamma_{\Sphere}\to\Gamma_{\Sphere}$, and we have then extended the graph map to a postsingularly finite topological exponential map (Lemma~\ref{lemma:extgraphmap}) which satisfies the four itemized properties in the theorem (Lemma~\ref{lem:embedding_satisfies}). Finally, Proposition~\ref{prop:thurston-uniqueness} shows that all topological exponential maps satisfying these conditions are Thurston equivalent.
\qed

\section{Holomorphic Exponential Maps}
\label{sec:topexp2holexp}
%Bottom left to bottom right arrow.

In the previous section, we have constructed for every preperiodic external address a topological exponential map satisfying the properties of Theorem~\ref{thm:extadd2topexpmap}. In this section, we are going to find a holomorphic exponential map with the same properties. It is here that we use the main result of \cite{HSS}:

\begin{thm}\capt{Characterization of Exponential Maps} \label{thm:HSS}%
A postsingularly finite topological exponential map is Thurston equivalent to a (necessarily unique) postsingularly finite holomorphic exponential map if and only if it does not admit a
Levy cycle.
\qedd
\end{thm}

Since the topological exponential map constructed in Theorem~\ref{thm:extadd2topexpmap} has no Levy cycle, we have now associated to every preperiodic external address a unique postsingularly finite holomorphic exponential map. In order to close the loop of the argument (see the diagram at the end of Section~\ref{sec:defns}), we need the following result.

\begin{prop}
\label{prop:closeloop} \capt{Dynamic Ray Lands at Singular Value}%
For every preperiodic external address $\s$, suppose that $\El$ is a holomorphic exponential map which is Thurston equivalent to the topological exponential map constructed in Theorem~\ref{thm:extadd2topexpmap} for external address $\s$. Then for $\El$, the dynamic ray at external address $\s$  lands at the singular value.
\end{prop}
\proof
We prove this claim by translating Thurston equivalence into the language of spiders, and using results from \cite{SZ2}.

For a postsingularly finite topological exponential map $f$, a {\em spider leg} is an injective curve $\gamma_i\colon[0,\infty]\to\Sphere$ with $\gamma_i(0)=e_i$ and $\gamma_i(\infty)=\infty$, where $e_1=0,e_2,\dots,e_{l+k}$ is the postsingular orbit. A {\em spider} is a collection of spider legs $\Phi=\{\gamma_1,\dots,\gamma_{l+k}\}$ which are disjoint except possibly for their endpoints. Two spiders are {\em equivalent} if there is an isotopy of $\Sphere$ relative to $P_f=\{e_1,\dots,e_{l+k},\infty\}$ which moves one spider to the other.

Every spider $\Phi$ lifts under $f$ to an image spider $\tilde\Phi=\{\tilde \gamma_1,\dots,\tilde\gamma_{l+k}\}$, where each $\tilde\gamma_i$ is the unique component of $f^{-1}(\gamma_{i+1})$ starting at $e_i$ (counting indices modulo the period as always). It is easy to check that equivalent spiders have equivalent image spiders, so the spider map acts on equivalence classes of spiders.

If $\gamma_1$ is a curve connecting $e_1=0$ to $\infty$ and is preperiodic under $f$ (as a curve) with preperiod $l$ and period $k$, then the $l+k$ legs on the orbit of $\gamma_1$ obviously form an invariant spider.

Now suppose that $f$ and $\El$ are Thurston equivalent: $\phi_1\circ f = \El\circ\phi_2$, where $\phi_1$ and $\phi_2$ are isotopic relative to $P_f$. Then $\Psi:=\phi_1\circ\Phi$ is a spider for $\El$. We can lift $\Psi$ under $\El$ to an image spider $\tilde\Psi$, and the condition $\phi_1\circ f = \El\circ\phi_2$ implies that $\tilde\Psi=\phi_2\circ\tilde\Phi$ (lifting spiders is compatible with Thurston equivalences). Since $\phi_1$ and $\phi_2$ are isotopic relative to $P_f$, it follows that the spiders $\Psi$ and $\tilde\Psi$ are equivalent: the map $\El$ has a fixed spider (up to equivalence) which is related by the Thurston equivalence to the preperiodic curves $\gamma_i$ of $f$.

Now \cite[Theorem~6.4]{SZ2} shows that this fixed spider of $\El$ can be replaced by an equivalent spider consisting only of dynamic rays (this theorem is part of the proof that every postsingularly finite exponential map has a dynamic ray landing at the singular value; the statement reads a bit differently in that context, but what the theorem actually does is to take a periodic spider given by \cite[Proposition~6.3]{SZ2} and turn it into a spider made of dynamic rays).

We now know that a dynamic ray landing at 0 has kneading sequence
$K(\s)=\It{\s}{\s}$, and the order of the image rays is as constructed in the
topological case. At this point, Algorithm \ref{algo:itinerary}
comes in handy: it allows us to reconstruct $\s$ uniquely from its kneading sequence $\u=\It{\s}{\s}$ and the order of the $\sigma^n(\s)$.

\begin{algo}
\label{algo:kneading&order->extaddr} \capt{Reversal of Construction of Kneading Sequence}
Suppose we are given the kneading sequence $\uu$ of a non-constant, non-periodic external address $\s$ and for all $n\ge 1$ we know the relative order of $\s$ and $\sigma^n(\s)$. Then the following algorithm uniquely recovers $\s$:
\begin{enumerate}
\item
For $n \geq 1$ define  $\delta_n := \left\{
\begin{array}{ll}
-1 & \mbox{ if } \sigma(\s) > \s \mbox{ and } \sigma^{\circ n}(\s) < \s \\
\ \ 1 & \mbox{ if } \sigma(\t) < \s \mbox{ and } \sigma^{\circ n}(\s) > \s \\
\ \ 0 & \mbox{ otherwise.}
\end{array} \right.$
\item Construct $\s$ as $s_n = \rm u_n - \delta_n$.
\end{enumerate}
\end{algo}

\proof
Since the relative order of $\sigma^{\circ n}(\s)$ and $\s$ is known for all $n$, we can compute $\delta_n$ for all $n$. We can thus reverse the computation in Algorithm~\ref{algo:itinerary}.
\qed

Therefore, for the holomorphic map $E_\lambda$ the dynamic ray at external address $\s$ lands at $0$. This concludes the proof of Proposition \ref{prop:closeloop}.
\qed

Now we can finish the proof of the first classification theorem.

\proofof{Theorem~\ref{thm:main}}
We have just finished the proof of the existence part of the theorem. For uniqueness, suppose there are two postsingularly finite exponential maps $\El$ and $\Elp$ which both have the property that their dynamic rays at external address $\s$ lands at the singular value. Then both maps have spiders consisting of this dynamic ray and its forward images, and there is a homeomorphism $\phi_1\colon\Cbar\to\Cbar$ which sends the spider of $\El$ to the spider of $\Elp$. Since $\phi_1(0)=0$ and both maps are exponential maps with asymptotic value $0$, it follows that $\phi_1$ lifts to another homeomorphism $\phi_2\colon\Cbar\to\Cbar$ so that $\phi_1\circ \El = \Elp\circ\phi_2$. The spiders assure that $\phi_1$ and $\phi_2$ are homotopic to each other relative to the postsingular set, so $\El$ and $\Elp$ are Thurston equivalent. By Theorem~\ref{thm:HSS}, it follows that $\lambda=\lambda'$.

Finally, the fact that every postsingularly finite exponential map is actually associated to a finite positive number of preperiodic external addresses is Theorem~\ref{thm:prepray@singval}.
\qed

A core ingredient in this proof was the combinatorial Algorithm~\ref{algo:kneading&order->extaddr} which allows to recover the external address $\s$ from its kneading sequence and the lexicographic order of the orbit of $\s$ under the shift. There are other ways to recover $\s$ which are perhaps more closely related to the dynamics of $f$; see for example the proof of \cite[Theorem~6.4]{SZ2}.

We have now proved Theorem~\ref{thm:main}: for every preperiodic
external address $\s$, there is a unique postsingularly finite
exponential map for which the dynamic ray at external address $\s$
lands at the singular value, and this describes all postsingularly
finite exponential maps.
\hide{Thus, association with the same $E_\lambda$
defines an equivalence relation on the space of preperiodic external
addresses. }
 It remains to determine combinatorially which external
addresses yield the same exponential map; this is
Theorem~\ref{thm:differentcodings}.

\proofof{Theorem \ref{thm:differentcodings}} Suppose $\s$ and $\s'$
give rise to the same exponential map $E_\lambda$: by definition,
the dynamic rays at addresses $\s$ and $\s'$ for $E_\lambda$ land
together at 0,  and hence they have the same itinerary with respect
to either dynamic partition, hence $\It {\s}{ \s'} = \It {\s'}{
\s'}$ and $\It{ \s'}{ \s} = \It{ \s }{\s}$. It also follows that
$\s$ and $\s'$ have the same period and the same preperiod.

Conversely, fix two preperiodic external addresses $\s$, $\s'$ of equal period and equal preperiod such that $\It{\s}{\s'} = \It{\s'}{\s'}$. Construct the holomorphic exponential map $E_\lambda$ corresponding to $\s'$ using Theorem \ref{thm:main}. $E_\lambda$ has a dynamic ray at external address $\s$, and since $\It {s}{ \s'} = \It {\s'}{ \s'}$, the rays at addresses $\s$ and $\s'$ land together at the same point $0$. By uniqueness, Theorem \ref{thm:main} constructs the same map from $\s$. Hence, both addresses correspond to the same holomorphic exponential map.

This shows the equivalence of Conditions (\ref{Item:Equiv_dynrays}), (\ref{Item:Equiv_itin1}), and (\ref{Item:Equiv_itin2}) in the claim. The equivalence of Conditions (\ref{Item:Equiv_dynrays}) and (\ref{Item:Equiv_pararays}) is Theorem~\ref{thm:landing_para_rays}. Before dealing with the remaining conditions, let us discuss the last statements of the theorem.

Dynamic rays landing at a common point always have equal period and
equal preperiod. The construction of kneading sequences makes it
clear that if $\s$ has preperiod $l$ and period $k$, then
$\It{\s}{\s}$ has preperiod $l$ and period $k'$ dividing $k$. Among
the external addresses $\s,\sigma(\s),\dots,\sigma^{k+l}(\s)$, only
the periodic ones can have equal itineraries with respect to $\s$,
and the number of those who do is obviously equal to $k/k'$.
Therefore, the corresponding dynamic rays land in groups of $k/k'$;
by \cite[Lemma~5.2]{S}, there can be no additional rays if $k>k'$,
while if $k=k'$, then the number of rays can be $1$ or $2$, and both
cases actually occur. Pulling back, this gives the number of
preperiodic dynamic rays landing at the singular value.

Condition (\ref{Item:Equiv_itin1&2}) implies the previous ones, so
we now show that if $\s$ and $\s'$ generate the same postsingularly
finite exponential map $\El$, then $\It{\s}{\s'}=\It{\s}{\s}$. For
every $j\ge 0$, the dynamic rays  $g_{\sigma^{j}(\s)}$ and
$g_{\sigma^{j}(\s')}$ land together and bound a component
$U_{j}\subset\C$ with real parts bounded below. If $\It{\s}{\s'}$
and $\It{\s}{\s}$ differ in their $k$-th entries, this implies that
$\El^{\circ k}(\ray)\subset U_0$ and, since both rays land together,
also $\El^{\circ k}(\rayp)\subset U_0$ and $U_k\subset U_0$. We show
that this is impossible.

There is a unique $k'\in\{0,1,\dots,k-1\}$ so that the first entry in $\sigma^{j}(\s)$ is the same as in $\sigma^{j}(\s')$, for every $j=k'+1,\dots,k-2$, but not for $j=k'$ (if there was no such $k'$, then the rays $\El^{\circ k}(\ray)$ and $\El^{\circ k}(\rayp)$ would be further apart than the rays $\ray$ and $\rayp$, and they could not be contained in $U_0$).

Then for $j=k'+1,\dots,k-2$, the restriction $\El\colon U_j\to U_{j+1}$ is a conformal isomorphism, and so is $\El\colon U_{k'}\to\C\sm\ovl U_{k'+1}$. It follows that $U_{k'+1}$ contains the singular value $0$ and hence the rays $\ray$ and $\rayp$, hence $U_{k'+1}\supset U_0\supset U_k$ (equality excluded). However, since the rays bounding $U_{k'}$ have more identical first entries in their external addresses than the rays bounding $U_k$, they must surround a smaller domain, and this is a contradiction.
\qed

%\addcontentsline{toc}{section}{Bibliography}

\end{document}